\documentclass[11pt,a4paper,twoside,openright,pdftex]{article}
\usepackage[english]{babel}
\usepackage{amsmath,amsthm,amssymb,amsfonts,epic,latexsym,hyperref,setspace}
\usepackage{natbib}
\usepackage[dvips]{graphicx}
\usepackage[utf8]{inputenc}
\usepackage{color}
\usepackage{xcolor}
\usepackage[margin=1in]{geometry}

\usepackage{tikz}
\usetikzlibrary{shapes,arrows,automata,mindmap,decorations}

\newcommand{\R}{\mathbb{R}}

\newcommand{\pr}{\mathbb{P}}

\newcommand{\esp}{\mathbb{E}}

\newcommand{\bgamma}{\boldsymbol \gamma}

\newcommand{\btheta}{\boldsymbol \theta}

\newcommand{\Y}{\mathcal Y}
\newcommand{\bx}{\textbf x}
\newcommand{\bZ}{\textbf Z}
\newcommand{\bY}{\textbf Y}
\newcommand{\mQ}{\mathbb{Q}}

\newcommand{\bb}{\boldsymbol \beta}

\newcommand{\bpi}{\boldsymbol \pi}

\newcommand{\I}{\mathcal{I}}

\newcommand{\calM}{\mathcal{M}}
\newcommand{\calP}{\mathcal{P}}

\graphicspath{{Figs/}}

\title{Modeling heterogeneity in random graphs through latent 
    space models: a selective review}
\author{Catherine Matias\footnote{Laboratoire de Math\'ematiques et Mod\'elisation d'\'Evry,
Universit\'e d'\'Evry Val d'Essonne,  UMR CNRS 8071, USC INRA, \'Evry,
France.  \texttt{email:}{catherine.matias@math.cnrs.fr}} 
\, \& \,
St\'ephane  Robin\footnote{AgroParisTech, UMR  518  Math\'ematiques et
  Informatique Appliqu\'ees, Paris 5\`eme, France.} \footnote{INRA, UMR 518 Math\'ematiques et Informatique Appliqu\'ees, Paris 5\`eme, France. \texttt{email:}{robin@agroparistech.fr}}}

\begin{document}
\DeclareGraphicsExtensions{.pdf, .jpg, .jpeg, .png}

\maketitle

\begin{abstract}
We present a selective review on probabilistic modeling of heterogeneity in random
graphs. We focus on latent space models and more particularly on
stochastic block models and their extensions that have undergone major developments in the
last five years.
\end{abstract}

\maketitle


\section{Introduction}

Network analysis arises in many fields of application such as biology,
sociology, ecology, industry, internet, etc. Random graphs represent a
natural way to describe how individuals or entities interact.  A 
network consists in a  graph where each node represents an
individual  and  an   edge  exists  between  two  nodes   if  the  two
corresponding individuals interact in  some way. Interaction may refer
to social  relationships, molecular  bindings, wired connexion  or web
hyperlinks,  depending  on  the  context.  Such  interactions  can  be
directed or not, binary (when only  the presence or absence of an edge
is recorded) or weighted (when a value is associated with each observed edge).

A huge literature exists on  random graphs and we refer the interested
reader e.g.   to the recent book by  \cite{Kolaczyk} for a general
and statistical approach to the field. A  survey  of  statistical
networks models appeared some  years ago in \cite{Goldenberg} and more
recently in \cite{Channarond_thesis,Snijders_review}.
In the present review, we focus on model-based methods for detecting
heterogeneity  in   random  graphs  and  clustering   the  nodes  into
homogeneous groups with respect to (w.r.t.) their connectivities.  This apparently specific focus still
covers a quite large literature.  The field  has  undergone so  many
developments in the past few years that there already exist other
interesting reviews of the field  and that the present one is supposed
to be complementary to those. In particular, we mention the complementary review
by \cite{Daudin_review_SFDS} focusing on binary graphs and the one by
\cite{JBL_SC} that reviews  methods  and  algorithms for
  detecting homogeneous subsets  of vertices  in binary or weighted and directed or undirected graphs. \\

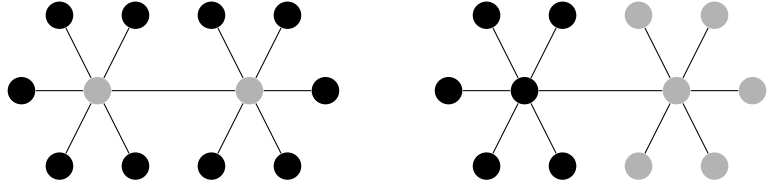
\begin{figure}[Ht]
  \centering
\begin{tikzpicture}
\tikzstyle{every node}=[circle,fill=white!70!black]
\node (1) at (0,0) {};
\node (2) at (2,0) {};
\tikzstyle{every node}=[circle,fill=black]
\node (S1) at (2.5,1) {};
\node (S2) at (2.5,-1) {};
\node (S3) at (1.5,1) {};
\node (S4) at (1.5,-1) {};
\node (S5) at (3,0) {};
\node (P1) at (0.5,1) {};
\node (P2) at (0.5,-1) {};
\node (P3) at (-0.5,1) {};
\node (P4) at (-0.5,-1) {};
\node (P5) at (-1,0) {}; 
\draw (1) to (2);
\draw (2) to (S1); \draw (2) to (S2);\draw (2) to (S3);\draw (2) to (S4); \draw (2) to (S5);   
\draw (1) to (P1); \draw (1) to (P2);\draw (1) to (P3);\draw (1) to (P4);  \draw (1) to (P5);
\tikzstyle{every node}=[]
\end{tikzpicture}
\hspace{1cm} 
\begin{tikzpicture}
\tikzstyle{every node}=[circle,fill=white!70!black]
\node (2) at (2,0) {};
\node (S1) at (2.5,1) {};
\node (S2) at (2.5,-1) {};
\node (S3) at (1.5,1) {};
\node (S4) at (1.5,-1) {};
\node (S5) at (3,0) {};
\tikzstyle{every node}=[circle,fill=black]
\node (1) at (0,0) {};
\node (P1) at (0.5,1) {};
\node (P2) at (0.5,-1) {};
\node (P3) at (-0.5,1) {};
\node (P4) at (-0.5,-1) {};
\node (P5) at (-1,0) {}; 
\draw (1) to (2);
\draw (2) to (S1); \draw (2) to (S2);\draw (2) to (S3);\draw (2) to (S4); \draw (2) to (S5);   
\draw (1) to (P1); \draw (1) to (P2);\draw (1) to (P3);\draw (1) to (P4);  \draw (1) to (P5);
\tikzstyle{every node}=[]
\end{tikzpicture}
  \caption{A toy example  of community detection (on the  right) and a
    more general  clustering approach based  on connectivity behavior
    (on the left) applied to the same graph with 2 groups. The two groups inferred are represented by colors (grey/black).}
  \label{fig:toy_cluster}
\end{figure}

The literature on  statistical approaches to random graphs  is born in
the social science community, where an important focus is given on
properties  such as  \emph{transitivity} ('the  friends of  your friends
  are likely to be your friends'), that is measured through clustering
indexes, and to \emph{community detection} that aims at detecting sets
of nodes  that share  a large number  of connections. However,  a more
general approach may be taken  to analyze networks and clusters can be
defined   as  a  set  of nodes  that  share  the  same
connectivity behavior. Thus,  we would like to stress  that there is a
fundamental  difference between general  nodes clustering  methods and
community  detection,  the  latter  being  a particular  case  of  the
former. 
For instance general nodes clustering methods  could put all the
hubs (i.e. highly connected nodes) together in one group and all peripheral nodes into another
group,  while  such  clusters  could never  be  obtained  from  a
community detection  approach. This is illustrated in  the toy example
from Figure~\ref{fig:toy_cluster}.

In the present  review, we do not discuss  general clustering methods,
but  only those  that  are  model-based. In  particular,  there is  an
abundant   literature  on   spectral  clustering   methods   or  other
modularities-based  methods for  clustering  networks  and we  only
mention those for which some performance results have been established
within heterogenous random graphs models. More precisely, a large part
of this review focuses on the stochastic block model (SBM) and we discuss
non  model-based  clustering  procedures  only when  these  have  been
studied in some way within SBM.

Exponential random graph models also constitutes a broad class of network models, which is very popular in the social science community (see \cite{GPK07} for an introduction). Such models have been designed to account for expected social
behaviors such as the {transitivity} mentioned earlier. {In some way,  they} induce some heterogeneity, but they
are not generative and the identification of some unobserved structure
is not their primary goal. For the sake of simplicity and homogeneity,
we decide not to discuss these models in this review.

The  review is  organized  as follows.  Section~\ref{sec:state_space} is  a
general  introduction to latent  space models  for random  graphs, that
covers {a large class of the  clustering methods that rely on a probabilistic model}.  
 More specifically, Section~\ref{sec:common} introduces the properties
shared       by       all       those      models     {while
  Sections~\ref{sec:latent_space} and~\ref{sec:other_latent} present
  models included in this class, except for SBM.} Then Section~\ref{sec:SBM}
focuses more precisely on the stochastic block model, which assumes that the
nodes are clustered into a finite number of different latent groups, controlling the
connectivities       between      the      nodes       from      those
groups.  {Section~\ref{sec:inference} is dedicated to
  parameter inference while Section~\ref{sec:clustering} focuses on clustering
  and model selection, both in SBM.}
In Section~\ref{sec:ext}, we review
the  various generalizations  that have  been proposed  to SBM  and in
Section~\ref{sec:conclusion}, we briefly conclude this review on discussing
the  next  challenges  w.r.t.  modeling heterogeneity  in  random
graphs {through latent space models}. 

Note that most of the following results are stated for undirected binary or
  weighted  random graphs with no  self-loops. However easy
  generalizations may  often be obtained for directed  graphs, with or
  without self-loops.

\section{Latent space models for random graphs}
\label{sec:state_space}
 \subsection{Common properties of latent space models}
\label{sec:common}
Let us first describe the set of observations at stake. 
When dealing with random graphs, we generally observe an adjacency matrix
$\{Y_{ij}\}_{1\le  i,j\le n}$ characterizing  the relations  between $n$
distinct individuals or nodes.  This can either be a binary matrix ($Y_{ij}\in \{0,1\}$ indicating presence or
absence of  each possible edge)  or a vector valued  table ($Y_{ij}\in
\R^s$ being  a -- possibly multivariate -- \emph{weight}  or any value  characterizing the
relation between nodes $i,j$).  The graph may be directed or undirected
(in which case $Y_{ij}=Y_{ji}$ for any $1\le i,j\le n$), it may either
admit self-loops or not ($Y_{ii}=0$ for any $1\le i \le n$).

{Latent} space models generally assume the existence of a latent random variable,
whose value characterizes the  distribution of the observation. In the
specific case of random  graphs, observations are the random variables
$Y_{ij}$    that    characterize    the   relation    between    nodes
$i,j$. {Note that
assuming that  the edges are independent variables  distributed from a
mixture model  - in which  case the latent variables  characterize the
edges behaviors - would not
make advantage of the  graph structure on the observations.} Thus, one
usually assumes that the latent variables rather characterize the nodes
behaviors  and each  observation $Y_{ij}$  will have  its distribution
characterized through the two different latent variables at nodes $i$ and $j$.
We thus assume that there exist some independent latent random variables
$\{Z_i\}_{1\le i \le n }$ being indexed by the set of nodes.  Moreover,
conditional on the $Z_i$'s,  the observations $\{Y_{ij}\}_{1\le i,j\le n}$
are independent and the distribution  of each $Y_{ij}$ only depends on
$Z_i$ and $Z_j$.  
This general  framework is considered  in \cite{BJR07}, where
  the  topological   properties  of  such  random   graphs  are  studied
  extensively  from a  probabilistic  point of  view.  Note that  this
  model results  in a set  of observations $\{Y_{ij}\}_{1\le
    i,j\le   n}$  that   are   not  independent   anymore.  In   fact,
  the dependency structure induced on the $Y$'s is rather complex as will be seen below.

Now,  we   will  distinguish  two  different  cases   occurring  in  the
literature: the latent random variables $Z_i$ may either take finitely
many values denoted by $\{1,\ldots, Q\}$, or being continuous and belong to some latent
space, e.g.  $\R^q$ or $[0,1]$.  In both cases (finite or
continuous), the network characteristics are summarized through a low dimensional latent space. 
The first case corresponds to the stochastic block model (SBM) and will
be   reviewed   in   extension  below   (see   Section~\ref{sec:SBM}).
The     second    case    will     be    dealt     with    in
  Sections~\ref{sec:latent_space} and~\ref{sec:other_latent}.\\

Before  closing this  section, we  would like  to explain  some issues
appearing when dealing with parameter
estimation  that are common  to all  these models.  Indeed, as  in any
{latent} space model, likelihood may
not be computed exactly except for small sample sizes, as it requires summing over
the set of possible latent configurations that is huge. 
Note that some authors  circumvent this computational problem  by considering the latent  variables   as   model  parameters   and   computing  a   likelihood
  conditional on these latent variables.  Then the observations 
  are independent and the corresponding likelihood has a very simple
  form  with  nice  statistical  properties.  As  a  counterpart,  the
  maximization  with respect  to  those latent  parameters raises  new
  issues. {For example, when the latent variables are discrete, this} results in a discrete optimization problem with associated combinatorial complexity. Besides
the  resulting  dependency  structure  on  the  observations  is  then
different {\citep[see for instance][]{Bickel_Chen,
    Rohe_Chat,Rohe_Yu}.}

\begin{figure}[Ht]
 \begin{tabular}{p{.33\textwidth}p{.33\textwidth}p{.33\textwidth}}
  \begin{tabular}{p{.4\textwidth}}
  $p(\{Z_i\}, \{Y_{ij}\})$ \\
   \includegraphics[width=.4\textwidth]{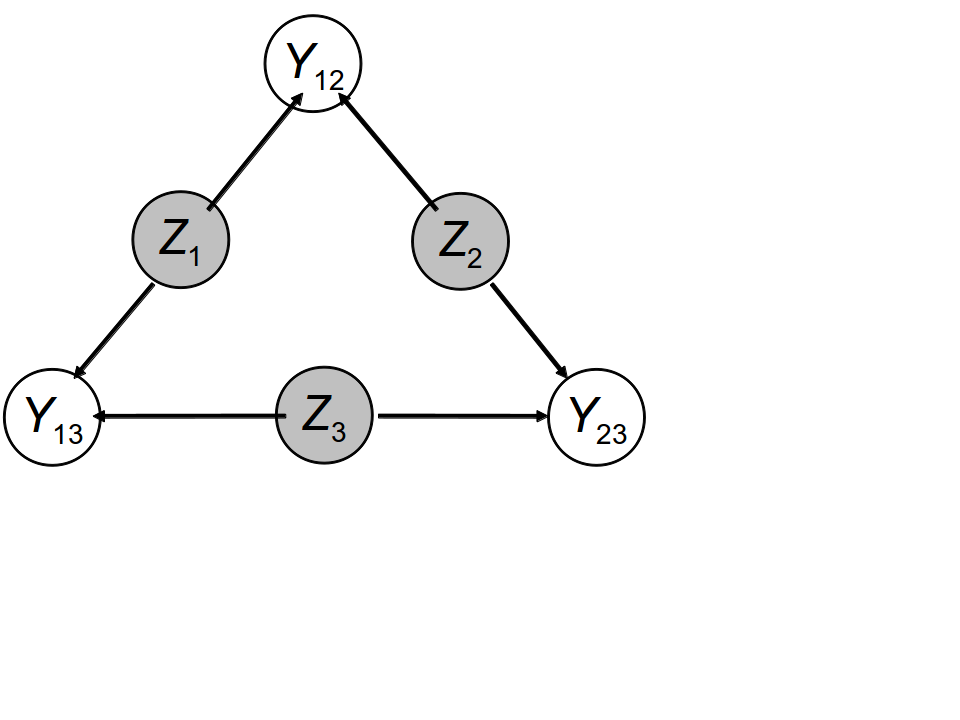} 
  \end{tabular}
  &
  \begin{tabular}{p{.4\textwidth}}
  Moralization of $p(\{Z_i\},  \{Y_{ij}\})$ \\
  \includegraphics[width=.4\textwidth]{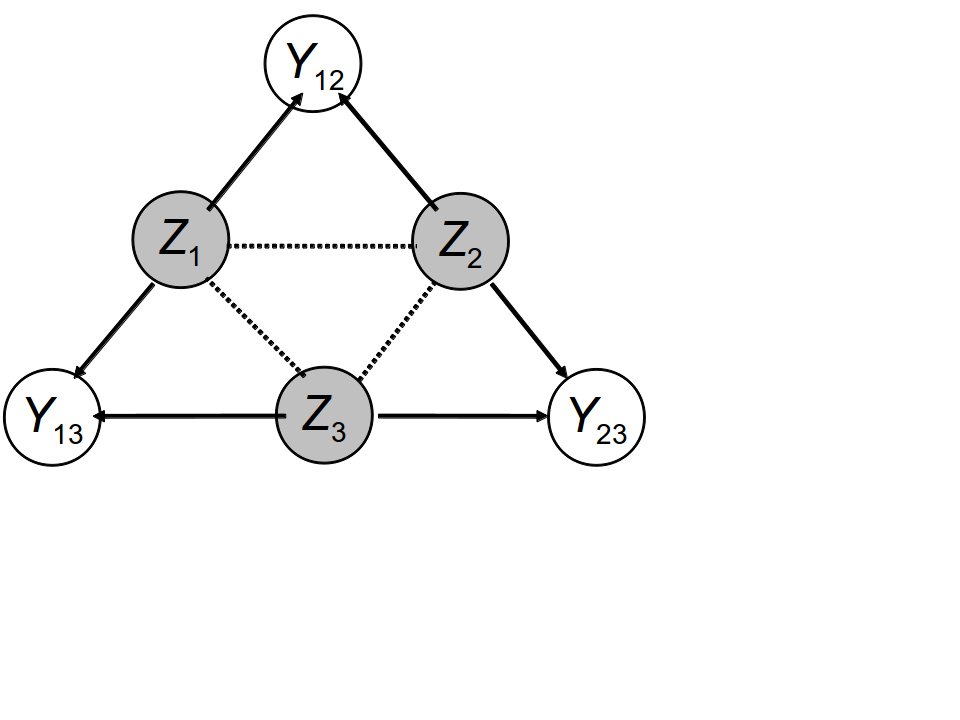}
  \end{tabular}
  &
  \begin{tabular}{p{.4\textwidth}}
  $p(\{Z_i\}| \{Y_{ij}\})$ \\
  \includegraphics[width=.4\textwidth]{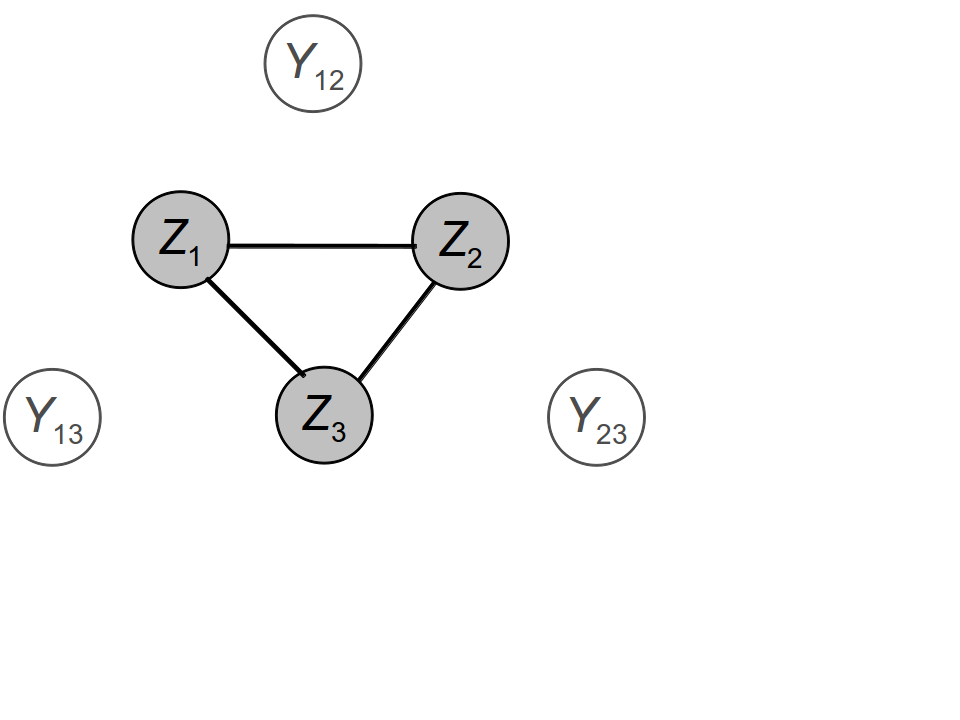}
  \end{tabular}
 \end{tabular}
 \caption{{Graphical  representation of  the  dependency structure  in
     {latent} space models for graphs.  Left:  {Latent} space model as a
     directed probabilistic  graphical model. Center:  Moralization of
     the  graph.   Right:  Conditional  distribution   of  the  latent
     variables    as    an    undirected    probabilistic    graphical
     model.  Legend: Observed variables (filled white), latent variables
     (filled gray) and conditioning variables (light gray lines).}
\label{fig:SBM}}
\end{figure}

However,   the   classical   answer  to   maximum   likelihood
  computation  with   latent  variables  lies   in  the  use   of  the
  expectation-maximization (\texttt{EM}) algorithm \citep{DLR}.
Though, the
\texttt{E}-step of the  \texttt{EM} algorithm  may be performed  only when
the  distribution of  the  latent variables  $Z_i$,  conditional on  the
observations $Y_{ij}$, may be easily computed.  This is the case for
instance in classical finite mixture  models (namely when observations are associated with respectively independent latent variables) as well as in models with more complex
dependency structure such as hidden Markov models \citep[HMM, see for
instance][]{Ephraim_Merhav,HMMbook} or more general conditional random
fields, where the distribution of  the latent $Z_i$'s conditional on the
observed         $Y_{ij}$'s          is         explicitely         modeled
\citep{Lafferty_etal,Sutton_CRF}. In  the case of  random graphs where
latent  random  variables  are  indexed  by the  set  of  nodes  while
observations  are indexed  by pairs  of nodes,  the  distribution  of  the $Z_i$'s  conditional  on  the
$Y_{ij}$'s   is  not  tractable.  
The reason  for this  complexity is
  explained in Figure~\ref{fig:SBM}. In this figure, the left panel reminds that the
  latent variables  $\{Z_i\}$ are  first drawn independently  and that
  the observed variables $\{Y_{ij}\}$ are then also drawn independently,
  conditional on the $\{Z_i\}$ with distribution that only depends on 
  their \emph{parent}  variables \citep[see][for  details on
    graphical models]{Lauritzen}.   The \emph{moralization step}  shows that
  the parents  are not independent anymore when  conditioning on their
  common offspring.  Indeed, if $p(Z_i,  Z_j, Y_{ij}) =  p(Z_i) p(Z_j)
  p(Y_{ij}|Z_i,  Z_j)$,  then  $p(Z_i,  Z_j  |Y_{ij})  =  p(Z_i,  Z_j,
  Y_{ij})/p(Y_{ij})$  can  not  be  factorized anymore  ('parents  get
  married').  The right  panel gives  the resulting  joint conditional
  distribution  of the $\{Z_i\}$  given the  $\{Y_{ij}\}$, which  is a
  clique.  This dependency  structure prevents  any  factorization, as
  opposed    to   models    such   as    HMM    or   {other graphical models} where this
  structure is tree-shaped.

Thus,
\texttt{EM} algorithm does  not apply and other strategies  need to be
developed for  parameter inference. These  may be classified  into two
main  types:  Monte  Carlo  Markov chains  (\texttt{MCMC})  strategies
\citep{MCMC}            and           variational           approaches
\citep{Jaakkola_tuto}.  The former methods aim at
sampling  from  the  true   conditional  distribution  of  the  latent
variables conditional on observed ones (e.g.
relying  on  a  Gibbs  sampler)  and  suffer  from  low  computational
efficiency. In fact, these methods are limited to network sizes of the
order  of a few  hundreds of  nodes. The
variational approaches result in \texttt{EM}-like algorithms in which an approximate conditional distribution of the latent variables conditional  on observed  ones is used. They are more efficient and may  handle larger data
sets \citep[up to few few thousands of nodes, 
e.g.][]{Zanghi_PatternRec08,Zanghi_AAS10}.  In
general, variational approaches suffer from a lack of convergence of the parameter
estimates {to the
  true parameter value, as the sample size increases} \citep{cv_varEM}. But in the
specific case of random graphs  {latent}  space models, they appear to be surprisingly
accurate. The  reason for this will  be given, at least  for SBM, in
Section~\ref{sec:clustering}  below.  The methods  for  parameter
inference in SBM will be described in
Section~\ref{sec:inference}.

\subsection{Latent space models (for binary graphs)}
\label{sec:latent_space}
Latent  space models  have been  developed  in the  context of  binary
graphs only. 
In   \cite{Hoff_etal}, the  latent
space $\R^q$ represents a \emph{social space} where the proximity of the actors 
induces a  higher probability of  connection in the graph.   Thus, only
relative positions in this latent space are relevant for the model. 
More  precisely, the  model is  defined for  binary random  graphs and
allows for  covariate vectors $\bx_{ij}$  on each relation  $(i,j)$. Two
different  parametrization have been  proposed in  \cite{Hoff_etal} to
deal with undirected and directed graphs, respectively. 
 For undirected  graphs, the  probability of connection  between nodes
 $i,j$ is parametrized through a logistic regression model 
\begin{equation*}
\text{logit}(\pr(Y_{ij}=1 |Z_i,Z_j, \bx_{ij})) =
  \frac{\pr(Y_{ij}=1 |Z_i,Z_j, \bx_{ij})}{1- \pr (Y_{ij}=1
    |Z_i,Z_j,\bx_{ij})} = \alpha +\beta^\intercal \bx_{ij} -\|Z_i-Z_j\|,
\end{equation*}
where $\|\cdot\|$ denotes Euclidean norm in latent space $\R^q$, 
the model parameters are $\alpha,\beta$ and $u^\intercal$ denotes the
transpose of vector $u$.  Note that the Euclidean norm could be
replaced by any kind of distance. 
In the case of directed networks, the distance is replaced by minus the
scalar product $Z_i^\intercal Z_j$, normalized by the length $\|Z_i\|$
of vector $Z_i$. Thus, the model becomes 
\begin{equation*}
\text{logit}(\pr(Y_{ij}=1 |Z_i,Z_j, \bx_{ij})) 
=     \alpha     +\beta^\intercal     \bx_{ij}
  +\frac{Z_i^\intercal Z_j}{\|Z_i\|} .
\end{equation*}
Note that in the distance case, the latent variables $\{Z_i\}$ might be recovered
only up to rotation, reflection and translation as these operations would
induce  \emph{equivalent   configurations}.   
Whether   these
  restrictions  are sufficient  for ensuring  the uniqueness  of these
  latent vectors has not been investigated to our knowledge.
Also note that in the model proposed here, the latent positions
  $Z_i$'s are considered as model parameters. Thus, the total number of
parameters  is $nq-q(q+1)/2  +2$  (including $\alpha$  and $\beta$), which
can  be  quite  large  unless  $q$ is  small.  {The  model  is
  specifically designed and thus mostly applied
  on social networks.}

\cite{Hoff_etal}  consider   a  Bayesian  setting   by  putting  prior
distributions  on   $\alpha,  \beta$  and   the  $Z_i$'s  and   rely  on
\texttt{MCMC} sampling to do parameter inference. The authors
  first   compute  a likelihood  that has  a  very simple  form (since  latent
  variables are considered as parameters, the observations are i.i.d.) and argue that
  this likelihood is convex w.r.t. the distances and may thus
  be first  optimized w.r.t.  these. Then, a  multidimensional scaling
  approach  enables  to identify  an  approximating  set of  positions
  $\{Z_i\}$ in  $\R^q$ fitting those distances.  These estimates $\hat
  Z_i$   form  an   initialization  for   the  second   part   of  the
  procedure. Indeed, in a second step, the authors 
  use an  acceptance-rejection algorithm to sample  from the posterior
  distribution of $(\alpha,\beta, \{Z_i\}_{1\le i \le n})$ conditional
  on the observations.
Note that with this latent space model, the nodes of the graph are not
automatically clustered into groups as  it is the case when the latent
space is  finite. For  this reason, \cite{Handcock_etal},  proposed to
further model the latent positions 
through a finite mixture  of multivariate Gaussian distributions, with
different   means    and   spherical   covariance    matrices.   Two
procedures are  proposed for  parameter inference: either  a two-stage
maximum likelihood method, where the first stage estimates the latent
positions   as   in  \cite{Hoff_etal}   (relying   on  a   simple-form
likelihood),  while the second  one is  an \texttt{EM}  procedure with
respect to the latent  clusters, conditionally on the estimated latent
positions; or a
Bayesian approach based on \texttt{MCMC} sampling. Besides, the number
of clusters  may be determined  by relying on  approximate conditional
Bayes factors.

The latent eigenmodel introduced in \cite{Hoff_Nips} defines the probability of
connection between nodes $(i,j)$  as a function of possible covariates
and a term of the form $Z_i^\intercal\Lambda Z_j$ where $Z_i\in \R^q$ is a
latent vector  associated to node $i$  and $\Lambda$ is  a $q\times q$
diagonal matrix. The author shows  that this form encompasses both SBM
and  the model  in  \cite{Hoff_etal}. However,  here  again, the  node
clustering is not induced by the model. The model is again applied to social
sciences but also linguistics (a network of word neighbors in a
text) and biology (protein-protein interaction network).

Latent  space models  have been  generalized into  \emph{random dot
  product graphs}.   Introduced     in     \cite{Nickel}     and
  \cite{Young_Scheinerman}, these models assume that each vertex is associated
  with a latent vector in $\R^q$ and the probability that two vertices
  are connected is then given by a function $g$ of the dot product of their respective
  latent  vectors. Three  different versions  of the  model  have been
  proposed in \cite{Nickel}, who shows that in at least two of those
  models, the resulting graphs obey a power law degree distribution, exhibit clustering, and
  have a low diameter.
In the model further studied in \cite{Young_Scheinerman}, each coordinate of those latent vectors $Z_i$ is drawn independently and identically from the
distribution  $q^{-1/2}  \mathcal{U}([0,1])^\alpha$, namely
for   any   $1\le   k   \le   q$,   the   coordinate   $Z_i(k)$   equals
$U_k^\alpha/\sqrt{q}$ where $U_1,\ldots, U_q$ are i.i.d.  with uniform 
distribution  $\mathcal{U}([0,1])$  on $[0,1]$  and  $\alpha>1$ is  some
fixed parameter. Moreover, the  probability of connection between two
nodes $i,j$ is exactly the dot product of corresponding latent vectors
$Z_i^\intercal Z_j$. 
 Interestingly, the  one-dimensional  ($q=1$)
version  of  this model  corresponds  to  a  graphon model  (see  next
section and Section~\ref{sec:graphon}) with function $g(u, v) = (uv)^\alpha$.
To our knowledge,  parameter inference has not been  dealt with in the
random dot product models (namely inferring $\alpha$ and in some models the parametric
link function $g$).  
However,   \cite{Tang_Sussman_Priebe}  have   proposed  a   method  to
consistently estimate
latent positions (up to an orthogonal transformation), relying on the eigen-decomposition of
$(AA^\intercal)^{1/2}$,  where  $A$ is  the  adjacency  matrix of  the
graph.  
We mention that the authors also provide classification results in a supervised setting where
latent positions are  labeled and a training set  is available (namely latent  positions  and their  labels  are observed).  However their  convergence results  concern a  model of  i.i.d. observations where     latent    variables     are     considered    as     fixed parameters.  {The  model has  been  applied  to  a graph  of web pages from Wikipedia in \cite{Sussman_Tang_Priebe}.}

Before closing this  section, we mention that the  problem of choosing
the dimension $q$  of the latent space has not been  the focus of much
attention.  We already  mentioned that  the number  of  parameters can
become quite large with $q$. In  practice, people seem to use $q=2$ or
$3$ and  heuristically compare the resulting fit  (taking into account
the number of parameters in each case). 
However the impact of the choice of $q$ has not been
investigated thoroughly. 
Moreover, as already mentioned, the parameters'   identifiability (with fixed $q$) has not been investigated in any of  the models described above.

\subsection{Other {latent} space models} \label{sec:otherstate}
\label{sec:other_latent}
 
Models   with   different  {latent}    spaces   have  also   been
  proposed.  In \cite{Daudin_etal10},  the latent  variables $\{Z_i\}$
  are supposed to belong to the simplex within $\R^Q$ (namely, $Z_i = (Z_{i1}, \dots, Z_{iQ})$ with all
$Z_{iq} > 0$ and $\sum_q Z_{iq} = 1$). 
As  for  the  inference,  the  latent positions  in  the  simplex  are
considered as  fixed and maximum  likelihood is used to  estimate both
the positions and the connection probabilities. Note that, because the
$Z_i$'s are  defined in a  continuous space, the  optimization problem
with respect to the $Z_i$'s is manageable. Conversely, as they are considered as parameters, the $Z_i$'s have to be accounted for in the penalized criterion to be used for the selection of the dimension $Q$.
This model can
be viewed as a continuous version of the stochastic block model (that will
be extensively discussed in the next section): in the stochastic block model
the  $Z_i$'s  would be  required  to belong  to  the  vertices of  the
simplex.  {We mention that it has been used in modeling
  hosts/parasites interactions \citep{Daudin_etal10}.}

A popular (and simple) graph model is the degree  sequence model in
  which  a  fixed  (i.e.  prescribed)  or  expected  degree  $d_i$  is
  associated to  each node. In  the fixed degree  model \citep{NSW01},
  the  random graph  is sampled  uniformly among  all graphs  with the
  prescribed  degree  sequence. In  the  expected degree  distribution
  model  \citep{PaN03,ChL02},  edges   are  drawn  independently  with
  respective  probabilities  $d_i d_j  /  \bar{d}$,  where $\bar{d}  =
  n^{-1} \sum_i d_i$. (Note that the degree-sequence must satisfy some
  constraints to ensure that $ d_i d_j/\bar{d}$ remains smaller than
  1.)  When the  $d_i$'s are  known, no  inference has  to  be carried
  out. When  they are not,  the expected degree model  may be
    viewed  as  a  latent  space  model  where  the  latent  variable
  associated to each node  is precisely $d_i$. {This model has
    been applied  to the structure of the  Internet \citep{PaN03}, the
    world wide web, as well  as collaboration graphs of scientists and
    Fortune 1000 company directors \citep{NSW01}.}

The graphon model (or $W$-graph) is a another popular model in the
probability community as it can be viewed as a limit for dense graphs \citep{Lovasz_Szegedy}. This model states that nodes are each associated with hidden variables $U_i$, all independent and uniformly distributed on $[0, 1]$. A graphon function $g: [0, 1]^2 \mapsto [0, 1]$ is further defined and the binary edges $(Y_{ij})_{1 \leq i < j \leq n}$ are then drawn independently conditional on the $U_i$'s, such that
\begin{equation} \label{eq:graphon}
 \pr(Y_{ij} = 1 | U_i, U_j) = g(U_i, U_j).
\end{equation}
The connections between this model and the stochastic block model will
be discussed in Section~\ref{sec:graphon}.


\section{Stochastic block model (binary or weighted graphs)}
\label{sec:SBM}
In this section, we {start  by  presenting  the
stochastic block model.} We
consider  a random  graph on  a set  $V=\{1,\ldots,n\}$ of  $n$ nodes,
defined as follows. Let $\bZ:= \{Z_1,\ldots, Z_n\}$ be i.i.d. random variables
taking  values  in  the  finite  set $\{1,\ldots,Q\}$  that  are  latent
(namely  unobserved),   with  some  distribution  $\bpi=(\pi_1,\ldots,
\pi_Q)$.   Alternatively,  each  $Z_i$  may be  viewed  as  a
  size-$Q$  vector  $Z_i=(Z_{i1},\ldots,   Z_{iQ})$  with  entries  in
  $\{0,1\}$ that sum up to one, whose distribution is multinomial $\mathcal{M}(1,\bpi)$.
The   observations   consist    in   the   set   $   \bY:=
\{Y_{ij}\}_{(i,j)\in \I}$ of random variables $Y_{ij}$ belonging to some space $\Y$, that characterize the
relations between  each pair of nodes  $i$ and $j$. When  the graph is
undirected  with no  self-loops, the  index set is  $\I=\{(i,j) ;  1\le i
<j\le n\}$, while it is equal to $\I=\{(i,j) ;  1\le i
\neq  j\le   n\}$  for  directed  graphs  with   no  self-loops.  Easy
generalizations   are  obtained   when  authorizing   self-loops.  The
distribution of the set of variables $\bY= \{Y_{ij}\}_{(i,j)\in \I}$ is as
follows: conditional on $\bZ= \{Z_i\}_{1\le i \le n}$, the $Y_{ij}$'s are
independent  and  the  distribution  of each  variable  $Y_{ij}$  only
depends on $Z_i$ and $Z_j$.  We let $F(\cdot; \gamma_{Z_iZ_j})$ denote
this distribution 
where $\bgamma =(\gamma_{q\ell})_{1\le q,\ell \le Q}$ is called the
connectivity  parameter.  Note  that  this matrix  is  symmetric  when
modeling undirected graphs. 
In conclusion, SBM is characterized by the following 
\begin{equation}
  \label{eq:SBM}
  \begin{array}{cl}
    \cdot  &\bZ  =  Z_1,  \ldots,  Z_n  \text{  i.i.d.  latent  random
      variables with distribution }
    \bpi \text{ on } \{1,\ldots,Q\}, \\
    \cdot &\bY  = \{Y_{ij}\}_{(i,j)\in \I} \text{  set of observations
      in } \Y^\I,\\
    \cdot & \pr (\bY | \bZ)    =\mathop{\otimes}_{(i,j)\in \I} \pr (Y_{ij} | Z_{i}, Z_j ) \text{ (conditional independence)}, \\
    \cdot & \forall {(i,j)\in \I} \text{ and } \forall 1\le q,\ell \le Q, \text{ we  have }  Y_{ij} |  Z_{i}=q, Z_{j}=\ell    \sim F(\cdot; \gamma_{q\ell}).
  \end{array}
\end{equation}
{SBMs  have  been  applied  to many  different  fields  and/or
  networks types, such as
  social sciences with social networks of individuals \citep{NS01,jernite2014}; 
  biology with regulatory transcription, metabolic, protein-protein interaction,
  foodwebs, cortex and hosts/parasites interaction networks
  \citep{Daudin_review_SFDS,Daudin_etal08,Picard_BMC,Mariadassou_10};  world wide web
  datasets    \citep{Latouche_overlap,Zanghi_AAS10};   cross-citations
  among journals \citep{Ambroise_Matias}; etc. 
}

Now  we  distinguish  binary  versions  of  the  model  from  weighted
ones.   Binary   SBMs   were   introduced  in   the   early   eighties
\citep{FH82,Holland_etal_83},  while weighted  versions  of the  model
appeared                 only                much                later
\citep{Mariadassou_10,Jiang_etal,Ambroise_Matias}. 
In  the  binary  SBM,  the  distribution of  $Y_{ij}$  conditional  on
$Z_i,Z_j$ is simply Bernoulli $\mathcal{B}(\gamma_{Z_iZ_j})$. Namely 
\begin{equation}
  \label{eq:SBM_binary}
\forall y \in \{0,1\}, \quad F(y; \gamma) =\gamma^y (1-\gamma)^{1-y}.
\end{equation}

Generalizing  the  model to  weighted  graphs,  we  consider that  the
distribution  of $Y_{ij}$  conditional  on $Z_i,Z_j$  is  any type  of
distribution  that depends only  on $Z_i,Z_j$.  More precisely,  it is
useful  to  restrict to  parametric  distributions,  such as  Poisson,
Gaussian, etc. {\citep[see for e.g.][]{Mariadassou_10}.} However, considering for instance an absolutely 
continuous distribution would induce a complete  graph, which is not desirable in
most applications. Thus, it makes  sense to consider instead a mixture
from a Dirac  mass at zero modeling absent  edges, with any parametric
distribution that models the strength or weight of present edges \citep{Ambroise_Matias}. For identifiability reasons, this latter distribution is restricted to have a cumulative distribution function (cdf) continuous at zero. In other words, we let 
\begin{equation}
\label{eq:SBM_weighted}
\forall y \in \Y, \quad F(y; \gamma) \sim \gamma^{1} G(\cdot, \gamma^{2}) +(1-\gamma^1) \delta_{0} (\cdot), 
\end{equation}
where  the  connectivity   parameter
$\gamma$ has  now two coordinates  $\gamma=(\gamma^1,\gamma^2)$ with $
\gamma^1 \in [0,1]$ and $G(\cdot,\gamma^2)$ is the
conditional distribution on the  weights (or intensity of connection),
constrained  to  have   a  continuous  cdf  at  zero.   When  all  the
$\gamma_{ql}^1$ are  equal to 1,  the graph is complete.  The particular
case  where  $G(\cdot,\gamma^2_{ql})$  is  the  Dirac mass  at  point  1
corresponds  to binary  SBM. Weighted  SBM may  for  instance consider
$G$ to be truncated Poisson, or a (multivariate) Gaussian, etc. Also
note  that in  the non  binary case,  the model  may be  simplified by
assuming $\gamma^1_{ql}$ is constant (equal to some fixed $\gamma^1$), so that
connectivity  is constant  throughout  the different  groups and  only
intensity varies.  

In what follows, the whole parameter of the SBM distribution is denoted by $\btheta=(\bpi,\bgamma)$.

Note that  a particular  case of SBM  is obtained when  considering an
\emph{affiliation} structure. Namely, the connectivity parameter $\bgamma$ takes only two different values, depending whether $q=\ell$ or not. In other words,
\begin{equation}
  \label{eq:affiliation}
  \forall 1\le q,\ell\le Q, \quad 
\gamma_{q\ell}=\left\{
  \begin{array}{cc}
   \gamma_{\text{in}} & \text{ when } q=\ell ,\\
   \gamma_{\text{out}} & \text{ when } q \neq \ell .
  \end{array}
\right.
\end{equation}
Note also that when considering a binary affiliation SBM and assuming moreover
that  $\gamma_{\text{in}}  \gg  \gamma_{\text{out}}$,  the  clustering
structure induced by the model corresponds exactly to clustering based
on the graph topology (namely  searching for sets of nodes that almost
form a clique).  That corresponds to community detection \citep{Fortunato}. 
As already mentioned, unconstrained SBM induces a node clustering that is much more general than community detection.  \\

To  conclude this  section, we  {mention} Szemer\'edi's
regularity Lemma
\citep{Szemeredi} as a potential motivation for SBM. 
Indeed, this lemma roughly states that every large enough graph
can be divided  into subsets of about the same size  so that the edges
between different  subsets behave almost  randomly {\citep[see
  e.g.][for a full presentation of the result]{On_Szeme}.}
{Note that the lemma
  is not of  a probabilistic nature.  Moreover, the  number of subsets
  whose existence is ensured by this result   may be arbitrarily large.}

\section{Parameter estimation in SBM}
\label{sec:inference}

\subsection{{Parameters' identifiability}}
We start this section  by discussing identifiability of the parameters
in  SBM.  As in  any  mixture  model, the  parameters  of  SBM may  be
recovered only up to a permutation on the groups labels. This is known
as \emph{identifiability up to label switching}. However, the issue of
whether  this  restriction  is   or  not  sufficient  to  ensure  the
parameters' identifiability has been overlooked in the literature for a
long time. In fact, the question was first solved only recently in the
particular case of  binary (undirected) SBM with only  $Q=2$ groups in
Theorem~7   from   \cite{ECJ}.   It   was  later   fully   solved   in
\cite{ident_mixnet} for (undirected) SBM,  both in binary and weighted
cases,    including   parametric   and    non-parametric   conditional
distributions  on the weights.  Note that  \cite{Celisse_etal} provide
another  identifiability  result valid  for  (directed or  undirected)
binary SBM. \\

\subsection{{Parameter estimation versus clustering}}
It is important to  note that at least two different approaches
may be considered when dealing with graphs under SBM. The first one is to estimate
the SBM parameters first, or at the same time as the nodes clusters. The second
one is  to cluster  the nodes  first (with no  information on  the SBM
parameters)  and  then  recover   the  SBM  parameters  through  these
estimated clusters. The latter is less related to SBM since generally,
the clustering  is done without  using the model  and is based  on the
graph structure.  In
Section~\ref{sec:clustering} , we discuss these 
methods when their theoretical properties within SBM have been discussed. As for
the  rest of  the  current section,  we  focus on  the first  approach
(namely parameter estimation in SBM). \\

\subsection{{\texttt{MCMC} approaches}}
\cite{SN97} have developed \texttt{MCMC} methods 
for  Bayesian  parameter  estimation  procedures in  binary  SBM.   We
mention that these
authors first considered a maximum likelihood method, that is limited to graphs
with up  to 20 or 30 vertices  for binary SBM with  $Q=2$ groups. Then
they proposed a Gibbs sampler for Bayesian estimation of the 
parameters.  
More   precisely,    given    current   values
  $(\bZ^{(t)}, \btheta^{(t)})$  of both the latent groups  and the SBM
  parameter, the algorithm samples 
  \begin{itemize}
  \item  $\btheta^{(t+1)}=(\bpi^{(t+1)},  \bgamma^{(t+1)})$  from  the
    posterior distribution given  the complete data $(\bZ^{(t)},\bY)$,
    namely by using 
    \begin{align*}
\pr(\bpi^{(t+1)}  =\bpi | \bZ^{(t)},  \bY) &  \propto \mu_{\bpi}(\bpi)
\prod_{i=1}^n \prod_{q=1}^Q \pi_{q}^{Z_{iq}^{(t)}} , \\
\pr(\bgamma^{(t+1)} =\bgamma| \bZ^{(t)}, \bY) & \propto \mu_{\bgamma}(\bgamma)
\prod_{(i,j)\in \I} \prod_{1\le q,\ell \le Q} F(Y_{ij} ; \gamma_{q\ell})^{Z_{iq}^{(t)}Z_{j\ell}^{(t)}},
    \end{align*}
where  $\propto$  stands  for  "proportional  to"  and  $\mu_{\bpi  },
\mu_{\bgamma}$  are the  prior distributions  on $\bpi$  and $\bgamma$
respectively;
\item  For $i=1$  to  $n$, sample  $Z_i^{(t+1)}$  from its  posterior
  distribution             given            $(\bY,Z_{1}^{(t+1)},\ldots
 , Z_{i-1}^{(t+1)},Z_{i+1}^{(t)}, $ $ \dots , Z_n^{(t)},
 \btheta^{(t+1)})$. To do this step, one uses the formula 
 \begin{equation}
   \label{eq:Gibbs}
\pr(Z_i=q| \bY,\{Z_{j}\}_{j\neq i }, \btheta) 
 \propto 
\pi_{q}\prod_{j ; (i,j)\in \I} \prod_{\ell =1}^Q F(Y_{ij} ; 
\gamma_{q\ell})^{Z_{j\ell}}.
 \end{equation}
  \end{itemize}

These results on Bayesian
estimation with a  Gibbs sampler have been extended  in \cite{NS01} to
handle directed graphs, an  arbitrary number of classes (restricted to
$Q=2$ earlier) and a  finite number of values for each relation
$Y_{ij}$. In practice, those  Bayesian methods are restricted to small
sample sizes (graphs with up to few hundred nodes). However,
very recent  attempts have been  made to develop  heuristic algorithms
with performances equivalent to exact \texttt{MCMC} procedures but much lower
running time \citep{Efficient_MCMC}. \\

\subsection{{Variational approximations}}
As  already  explained  in Section~\ref{sec:common},  the  \texttt{EM}
algorithm may not be performed exactly in SBM due to the intricate form of the
conditional  distribution of the  groups given  the data.  The natural
solution in such a case is to replace this conditional distribution by
its best approximation within a reduced class of distributions with simpler form.  
This leads to what is called a variational approximation  to the maximum
likelihood computation.  Let us explain this in more details. 
The data
log-likelihood $\mathcal{L}_{\bY}(\btheta)$ may be decomposed as 
follows
\begin{equation*}
  \mathcal{L}_{\bY}(\btheta) 
  :=               \log
  \pr(\bY; \btheta) 
= \log \pr(\bY, \bZ; \btheta) - \log \pr( \bZ| \bY ;\btheta)
\end{equation*}
and by taking on both sides  the expectation with respect  to some
distribution $\mQ$ acting only on $\bZ$, we get 
\begin{equation}\label{eq:decomp_loglik_esp}
    \mathcal{L}_{\bY}(\btheta)        =        \esp_{\mQ}(        \log
    \pr(\bY, \bZ ; \btheta) ) 
+\mathcal{H}( \mQ)
+ \mathcal{KL}( \mQ \|  \pr(  \bZ|  \bY;
\btheta) ) , 
\end{equation}
where     $\mathcal{H}(\pr)$     is     the     entropy     of     
distribution $\pr$ and $\mathcal{KL} (\pr \| \mQ)$ is the Kullback-Leibler
divergence between distributions $\pr$ and $\mQ$.
Starting from this relation, \texttt{EM} algorithm is an iterative
procedure based on the iteration  of the two following steps. Starting
from current parameter value $\theta^{(t)}$, we do 
\begin{itemize}
\item \texttt{E}-step: maximise the quantity $ \esp_{\mQ}(        \log
    \pr(\bY, \bZ ; \btheta^{(t)}) ) +\mathcal{H}(\mQ )$
    with respect to  $\mQ$. From~\eqref{eq:decomp_loglik_esp}, since $
    \mathcal{L}_{\bY}(\btheta^{(t)})  $ does  not depend  on $\mQ$,  this is
    equivalent to minimizing $\mathcal{KL}( \mQ \|  \pr(  \bZ|  \bY;
\btheta^{(t)}) )$ with respect to $\mQ$. The optimal solution is thus given by 
the  conditional   distribution  $  \pr(  \bZ|  \bY;  \btheta^{(t)})$  for current  parameter  value
  $\btheta^{(t)}$; 
\item \texttt{M}-step: keeping now  $\mQ$ fixed, maximize the quantity
  $ \esp_{\mQ}(  \log \pr(\bY, \bZ ; \btheta)  ) +\mathcal{H}(\mQ)$
with respect to $\btheta$ and update the parameter value $\btheta^{(t+1)}$ to this
    maximiser.  
   As $\mQ$ does not  involve the parameter $\btheta$ , this
      is   equivalent  to   maximizing  the   conditional  expectation
      $\esp_\mQ( \log  \pr(\bY, \bZ ; \btheta))  $ w.r.t.  $\btheta$.
Note that  here, with our choice of  $\mQ$, this quantity
      is the      conditional expectation $\esp( \log \pr(\bY, \bZ ; \btheta) |\bY, \theta^{(t)}) $ w.r.t. $\btheta$.    Moreover, this will
    automatically increase the log-likelihood $ \mathcal{L}_{\bY}(\btheta)$.
\end{itemize}
When the
true distribution  $\pr(\bZ|\bY)$ is intractable  (e.g. when it
can not be factorized in any way), the exact solution from \texttt{E}-step may not be
computed.  Instead,  going  back  to~\eqref{eq:decomp_loglik_esp}  and
using that the Kullback-Leibler divergence term is positive, we obtain the
following lower bound 
\begin{equation}
  \label{eq:lower_bound}
      \mathcal{L}_{\bY}(\btheta) \ge     \esp_{\mQ}(        \log
    \pr(\bY,   \bZ  ;
    \btheta) ) + \mathcal{H}(\mQ). 
\end{equation}
In  the  variational approximation,  instead  of  computing the  exact
solution at \texttt{E}-step, we  rather search for an optimal solution
within  a restricted  class of  distributions, e.g.  within the
class of factorized distributions 
\[
\mQ(\bZ) = \prod_{i=1}^n \mQ(Z_i), 
\]
and the \texttt{M}-step is unchanged, with $\mQ$ fixed to the previous
solution.  
In the case of SBM, taking $\mQ$ within the class of factorized distributions:
\[
\mQ(\bZ)  =   \prod_{i=1}^n  \mQ(Z_i)   =  \prod_{i=1}^n
\prod_{q=1} ^Q \tau_{iq}^{Z_{iq}},
\]
where $\tau_{iq} = \mQ(Z_i=q)$ (with $\sum_q \tau_{iq} = 1$ for all $i$), the solution to the \texttt{E}-step at
the  current parameter  value  $\btheta$, within  the  above class  distributions satisfies the following fixed point relation (see
Proposition 5  in~\cite{Daudin_etal08} in the binary  case and Section
4.2 in~\cite{Mariadassou_10} for weighted graphs):
\[
\tau_{iq}   \propto   \pi_q  \prod_{   {j;  (i,j)\in   \I}}
\prod_{\ell=1}^Q [f(Y_{ij}; \gamma_{q\ell})]^{\tau_{j\ell}} .
\]
The  resulting approximation  is sometimes  called a  {\em  mean field
  approximation}   because,   when   considering   the   (conditional)
distribution of $Z_i$,  all other $Z_j$'s are set  to their respective
(conditional) means  $\tau_{iq}$. A link  can be made with  {\tt MCMC}
techniques  applied  to  SBM,  which   most  often  rely  on  a  Gibbs
sampler.  In  this Gibbs  sampling  step,  the  $Z_i$ are  iteratively
sampled   conditionally  on  the   observed  variables,   the  current
parameters and  the other $Z_j$'s, {that  is with distribution
  given by~\eqref{eq:Gibbs}. }
In the variational framework, the probability $\tau_{iq}$ is the best possible approximation of $\pr(Z_i=q  | \bY,  \btheta)$ in terms of Kullback-Leibler divergence, within the set of factorized distributions.

By doing this variational approximation, we only optimize the lower bound on the right hand-side
of~\eqref{eq:lower_bound}  with  respect   to  $\theta$  and  have  no
guarantee of approximating the  maximum likelihood estimator. In fact,
as  already mentioned  in Section~\ref{sec:common},  these variational
approximations  are  known  to  be  non convergent  in  regular  cases
\citep{cv_varEM}.  More  precisely,  in  general  the  \texttt{E}-step
approximation  prevents   convergence  to   a  local  maxima   of  the
likelihood. 
However, it appears that for SBM, both empirical and
theoretical results  ensure that these procedures  are pretty accurate
(see Section~\ref{sec:asympto} below as well as
Section~\ref{sec:clustering} for more details).

In  \cite{Daudin_etal08},   a  frequentist  variational   approach  is
developed in the  context of binary SBM, while  the method is extended
to weighted graphs in~\cite{Mariadassou_10}. In~\cite{Picard_BMC}, the
variational  procedure is  applied in  the  context of  binary SBM  to
different biological networks, such as a transcriptional regulatory
network, a metabolic network, a cortex network, etc. 
A  Bayesian  version  of  the  variational  approximation  appears  in
\cite{Latouche_Bayes_EM} for  binary graphs \citep[see  also][for more
details]{Latouche_thesis} and more recently 
  in \cite{Aicher_etal} for weighted ones.  The aim is now to approximate the posterior
distribution $\pr(\bZ,\btheta | \bY)$ of both the groups and the
parameters,  given the  data,  a  task realized  within  the class  of
factorized distributions.  
Online extensions of the  (frequentist) variational approach have been
given  in \cite{Zanghi_PatternRec08,Zanghi_AAS10}  for binary  SBM and
more general weighted graphs where the conditional distribution of the
edges, given the nodes groups belongs to the exponential family.

Note that  efficient implementations  of the variational  approach for
binary  or weighted  graphs  are available  in  several softwares  and
packages:  \texttt{MixNet}   \citep{Picard_BMC}  for  binary   graphs,  \texttt{WMixnet}
\citep{JBL_soft}  for  both  binary  and  weighted  graphs,  \texttt{Mixer}
\citep[R package,][]{Mixer}  for binary  SBM and \texttt{OSBM}  for overlapping
extensions \citep[R package,][]{OSBM}. 
\\

\subsection{{Moment methods}} 
The very first work on stochastic block models is due to
\cite{FH82} where parameter estimation is  based on a moment method. The
article  considers  only  a  binary  affiliation  model  and  aims  at
estimating      the      parameters      $\gamma_{\text{in}}$      and
$\gamma_{\text{out}}$  under  the  assumption that  group  proportions
$\bpi$ are known. Note that this method has some drawbacks and has been both
discussed      and     extended      in      Section~3.1      from
\cite{Ambroise_Matias}.  
A  similar  approach  based  on  moment  methods  has  been  taken  in
\cite{BCL11}  for  a more  general  model,  namely  the graphon  model
described below in Section~\ref{sec:graphon}.\\

\subsection{{Pseudo or composite likelihood methods}}
Besides  considering  moment methods  for  binary affiliation  graphs,
\cite{Ambroise_Matias} propose {two different composite likelihood
approaches that are suited for binary and weighted affiliation graphs,
respectively.  In the weighted affiliation case,} their approach relies on optimizing a criteria that would    correspond to the
  the log-likelihood of  the $Y_{ij}$'s in a model where
  these  variables were  independent.   {As for  the binary  affiliation
  case, considering the same criteria  is not possible because in such
  a model the $Y_{ij}$'s are  a mixture of two Bernoulli distributions
  (one for out-group connections and the other for in-group ones),
  whose  parameters  cannot  be  identified. However,  by  looking  at
  triplets   $(Y_{ij},Y_{jk},Y_{ik})$   and   considering   these   as
  independent,  they  obtain  a  multivariate  Bernoulli  distribution
  whose parameters can be identified \citep{ECJ}.}
They  prove convergence  results
  justifying this pseudo-likelihood approach (see below) and exhibit a
  good accuracy 
  on simulations. In the same way, \cite{ACB13}
  proposed  a pseudo-likelihood approach  for binary  graphs. Starting
  with an initial configuration (namely nodes groups assignment), they
  consider the random variables $\mathbf{b}_i=(b_{iq})_{1\le q \le Q}$, where
  $b_{iq}$ is the number of connections of
  node $i$ to  nodes within class $q$ for the  given assignment of nodes
  classes.  They consider a pseudo (or composite) likelihood of these 
  variables   (namely   doing   as   if  the   $\mathbf{b}_i$'s   were
  independent)    and     optimize    this    criterion     with    an
  \texttt{EM}-algorithm.  At  the end  of  this  \texttt{EM} run,  the
  resulting clustering  on the nodes is  used as a  starting point for
  the next \texttt{EM} run.  Thus, the \texttt{EM} procedure is itself
  iterated many times until convergence. 
The  method  is   shown  to  perform  well  on   simulations  (but  no
theoretical result supports this approach w.r.t. parameter estimation).\\

\subsection{{Ad-hoc methods: nodes degrees}}
In \cite{Channarond}, a method based on the degree distribution of the
nodes  is explored for  (undirected) binary  SBM.  More  precisely, by
letting   $\bar  \gamma_q=\sum_{l=1}^Q   \pi_l  \gamma_{ql}$   be  the
probability of  connection of a node  given that it  belongs to class
$q$, the separability assumption  ensures that all the $\bar \gamma_q,
1\le q  \le Q$ are distinct. Under  this assumption, \citeauthor{Channarond}
propose a method  to estimate the parameters, based  only on the nodes
degrees. The approach is a generalization of a proposal by \cite{SN97}
in   the   case   of   $Q=2$    nodes   (see   Section   5   in   that
reference). Theoretical  results associated  with this method  will be
discussed in Section~\ref{sec:clustering}. Note that since the method only relies on the nodes
degrees, it is very fast and may easily handle very large graphs.
 \\

\subsection{{Asymptotic properties for parameter estimates}}
\label{sec:asympto}  
Very few convergence
results  have   been  obtained  concerning   the  different  parameter
estimation  procedures developed  in SBM.  Concerning  the variational
estimates, as previously said these estimators were not expected to be
convergent,  in  the sense  that  infinitely  many  iterations of  the
variational algorithm would not
necessarily  lead  to  a  local  maximum of  the  likelihood.  However
empirical results exhibit the accuracy of these variational
estimates \citep{Gazal}.   This latter  reference studies in  fact the
empirical convergence (from an  asymptotic perspective, as the size of
the graph increases) of  three different procedures: the (frequentist)
variational   estimate,  a   belief  propagation   estimate   and  the
variational Bayes  estimator.  The nice convergence  properties of the
variational estimator might be the joint consequence of
two different facts: 
\begin{enumerate}
\item the variational procedure  approximates a local maxima of
  the likelihood, as the number of iterations increases;
\item the maximum likelihood estimator is convergent to the
  true parameter value, as the size of the graph increases.
\end{enumerate}
Point (1) is  in fact a consequence of some kind  of degeneracy of the
model, where asymptotically the conditional distribution of the latent
variables $\bZ$ given the observed ones $\bY$ is a product of Dirac masses,
thus a factorized distribution and the variational approximation turns
out to  be  asymptotically  exact.   We  discuss  this  point  further  in
Section~\ref{sec:clustering}. \\
Now, point (2) has  been established in \cite{Celisse_etal} under some
assumptions. More  precisely, for binary (possibly  directed) SBM, the
authors prove  that maximum  likelihood and variational  estimators of
the group connectivities $\bgamma$ are consistent (as the size of the
graph increases).  However, they can not establish  the convergence of
the  same estimators  for the  groups proportions  $\bpi$  without the
additional  assumption that  the estimators  of $\bgamma$  converge at
rate faster than $\sqrt{\log n}/n$ (where $n$ is the number of nodes).
Note  that such  an  assumption is  not  harmless since  the rates  of
convergence of those estimates are still unknown.  Moreover, up to the
logarithmic   term,  the  rate   required  here   is  $1/n$   and  not
$1/\sqrt{n}$, which would
correspond to the parametric rate for an amount of $n^2$ (independent)
data. The fact that the
group  proportion parameters  $\bpi$ and  the  connectivity parameters
$\bgamma$ fundamentally  play a different  role in SBM occurs  in many
problems,  as  for  instance   for  the  model  selection  issue  (see
Section~\ref{sec:selection}). \\
In  the affiliation  SBM, \cite{Ambroise_Matias}  obtained convergence
results for the moment estimators  they proposed in the binary case as
well as  for the maximum composite likelihood  estimators developed in
both binary and  weighted cases. Note that  in this reference, the authors establish
rates  of  convergence  of  the procedures.  Surprisingly,  the  rates
obtained  there   are  not  of   the  order  $1/n$   but  $1/\sqrt{n}$
instead.  More precisely, in  the more general affiliation cases, they
establish an
asymptotic  normality result with  non degenerate  limiting covariance
matrix ensuring  that the estimators converge at  the usual $1/\sqrt{n}$
parametric rate (for an amount of $n$ "independent" data).  Then, in very specific subcases (e.g. equal
group proportions),  the limiting  covariance degenerates and  rate of
convergence increases to $1/n$. 
The issue  of whether these  rates are or  not optimal in  this context
remains open.

{We  conclude this section  by mentioning  that all  the above
  asymptotic results are 
  established only in a
  dense regime where the number of edges in the graph increases as $n$ grows to
  infinity.   Other   setups   such   as   letting   $Q$   fixed   but
  $\bgamma=\bgamma_n$  goes  to zero,  or  letting  $Q=Q_n$ increase  to
  infinity and $\bgamma$ fixed need to be further investigated.}
In  the  next  section,  we  discuss  asymptotic  properties  of  some
clustering  procedures.  Note   that  procedures  that  asymptotically
correctly recover the nodes groups (e.g. with large probability, w.l.p.) and base
their   parameter   estimation   on   these  estimated   groups   will
automatically be consistent (e.g.  w.l.p.). We refer to Theorem
4.1 in \cite{Channarond} for a formal proof of such a result.

\section{Clustering and model selection in SBM}
\label{sec:clustering}

Before starting this section, let us recall that clustering within SBM
is not limited to community detection. The latter corresponds to the very special case of an
  affiliation  model,  with  additional  constraint  that  intra-group
  connectivity $\gamma_{\text{in}}$ should  be larger than outer-group
  connectivity $\gamma_{\text{out}}$. Many  methods have been proposed
  to cluster the nodes within SBM, among which we distinguish maximum a
  posteriori   (MAP)  estimation   based  on   the   groups  posterior
  distribution given the data, from other methods. In the 
  first approach, parameters are estimated  first (or at the same time
  as the clusters)  and one considers the properties  of the resulting
  posterior distribution $\pr(\bZ|\bY ; \hat \theta)$ at an
  estimated  parameter  value  $\hat  \theta$,  while  in  the  second
  approach, the clusters are estimated first (without relying on parameter
  inference)  and then  parameters estimators  are  naturally obtained
  from these clusters. \\

\subsection{{Maximum a posteriori}} 
As previously  mentioned, \cite{Celisse_etal} studied  the behavior of
the maximum likelihood and the variational estimators in (binary) SBM. To this
aim, they have studied the posterior distribution of the groups, given
the data. These authors establish two different results. The first one
\citep[Theorem  3.1  in][]{Celisse_etal}   states  that  at  the  true
parameter  value, the  groups  posterior distribution  converges to  a
Dirac mass  at the actual  value of groups  configuration (controlling
also the corresponding rate of convergence). This result is valid only
at the true parameter value and not an estimated one.
The  second  result they  obtain  on  the  convergence of  the  groups
posterior  distribution \citep[Proposition 3.8  in][]{Celisse_etal} is
valid  at  an  estimated  parameter  value,  provided  this  estimator
converges at  rate at least $n^{-1}$  to the true value.  Note that we
already  discussed rates  of convergence  for SBM  parameters  and the
latter property has not been proved for any estimator yet. The article
\cite{Maria_Matias}  is more  dedicated  to the  study  of the  groups
posterior distribution in any  binary or weighted graph (their results
being in fact valid for  the more general latent block model described
in Section~\ref{sec:bipartite}).  The authors study this posterior for any
parameter value in the neighborhood  of the true value, thus requiring
only consistency  of a parameter estimator.  They establish sufficient
conditions for  the groups posterior distribution to  converge (as the
size of  the data  increases) to  a Dirac mass  located at  the actual
(random)  groups   configuration.   These  conditions   highlight  the
existence  of particular  cases  in SBM,  where some  \emph{equivalent
  configurations} exist,  and exact recovery  of the latent  groups is
not possible. 
They  also give  results  in a  sparse
regime  when the  proportion of  non-null entries  in the  data matrix
converges to zero. 

Note  that  those  results  of  convergence of  the  groups  posterior
distribution  to  a (product  of)  Dirac  mass(es)  at the  actual  groups
configurations  explains  the  accuracy  of  the  variational
approximation.  Indeed, {as the sample size increases, the error made by
  the variational approximation tends to zero}.\\

\subsection{{Other clustering methods}} 
Many other methods have been used to cluster  a graph and we
only consider those whose properties have been studied under
SBM  assumption.   Among  those  methods,  there  are  some  based  on
modularities, some based on the nodes degrees and we also mention some
properties of spectral clustering in a framework related to SBM.

In \cite{Bickel_Chen}, the authors show that groups estimates based on the use of  different  modularities  are  consistent  in  the  sense  that  with
  probability   tending   to   one,   these   recover   the   original
  groups of a binary SBM.  Quoting \cite{Bickel_Chen}, \emph{the Newman-Girvan
    modularity  measures  the fraction  of  edges  on  the graph  that
    connect vertices  of the  same type (i.e.  within-community edges)
    minus the expected  value of the  same quantity on  a graph
    with the  same community divisions but  random connections between
    the vertices}.  This modularity  is clearly designed for community
  detection purposes. 
The  authors  also introduce  a  \emph{likelihood
    modularity}, that is a  profile likelihood, where the nodes groups
  are considered  as parameters.  Under a  condition {on these
    modularities} that is quite
  difficult to  understand and whose consequences  remain unclear (see
  Condition I  in that reference),  they establish consistency  of the
  clustering   procedures  that   rely  on   these   modularities.  In
  {the} particular {case of Newman-Girvan
  modularity}, it  is  unlikely  that  Condition   I  is  satisfied in a non affiliation SBM.  
Moreover, the  likelihood modularity (or  profile likelihood)
  is computed  through a  stochastic search over  the node  labels. In
  practice, this  might raise  some issues that  are not  discussed by
  the authors.
\\

In the specific case of a binary graph with $Q=2$ nodes groups, \cite{SN97} proved that
the groups could be recovered exactly with probability tending to 1,
by using  a method based on the  nodes degrees (see Section  5 in that
reference). 
\cite{Channarond} generalized  the method to binary graphs
with any number of  groups.  Under the separability assumption already
mentioned above, the authors establish a bound on the probability of
misclassification of at least one node. \\

\cite{Rohe_Chat} propose a classification algorithm based on spectral
clustering that  achieves vanishing classification error  rate under a
model called binary SBM. In fact,  it is worth noticing that the setup
is the one of \emph{independent} Bernoulli random observations: latent
groups are viewed as parameters instead  of random variables.
This is strictly different from SBM.  
 In  the same  vein, \cite{Rohe_Yu}  are concerned  with a
framework in  which nodes of  a binary graph  belong to two  groups: a
\emph{receiver} group and a \emph{sender} group.  This is a refinement
of standard SBM, which assumes equal sender and receiver groups, and is motivated by the
study   of   directed  graphs.   They   generalize   the  results   of
\cite{Rohe_Chat}    to    a    framework    called    \emph{stochastic
  co-block model}. Here  again, this would be a  generalization of SBM,
except that they consider that the edges are independent random variables. 
The results from \cite{Rohe_Chat,Rohe_Yu} allow the number of groups to
grow with network size (i.e. nodes number) but require that node degrees increase
nearly linearly with this size, an assumption that can be
restrictive.   Note  that   \cite{CWA12}   provided  results   for
likelihood-based clustering in  a sparser  setup where  nodes degrees
increase poly-logarithmically w.r.t. the number of nodes. But here again,
the setup is the one of independent Bernoulli random variables and not
SBM.  Also  note  that  the  complexity  of  the  spectral  clustering
algorithm is  $O(n^3)$, which makes its  practical use computationally
demanding for large networks,  even though faster approximate versions
exist. \\

\subsection{\textbf{Model selection}}
\label{sec:selection}

As in most discrete  {latent} space models, the number of classes $Q$ is unknown in general and needs to be estimated. Few model selection criteria have been proposed up to now to address this question. \cite{NS01} introduced an information measure and a posterior entropy parameter that both can be used to evaluate the reliability of the clustering, but with no explicit model selection procedure as for the choice of $Q$. \\
For undirected binary graphs, \cite{Daudin_etal08} derived an ICL-like  (integrated complete likelihood) criterion. The ICL criterion has been first proposed by  \cite{BCG00} in the context of mixture models and is the same as the BIC (Bayesian information criterion) with an additional penalty term which corresponds to the entropy of the conditional distribution $\pr(\bZ|\bY)$.
\citeauthor{Daudin_etal08} used the term $\esp_\mQ(\log \pr(\bY, \bZ; \btheta))$ from the lower bound \eqref{eq:lower_bound} as a proxy for the expectation of the complete log-likelihood. Interestingly, they end up with a two-term penalty with the form
\[
\frac12 \left((Q-1) \log n + \frac{Q(Q+1)}{2} \log \frac{n(n-1)}{2} \right) , 
\]
where the first term refers to $\bpi$ and the second to $\bgamma$. This form reminds that the relevant sample size for the estimation of the group proportions $\pi_q$ is the number of nodes, whereas it is the number of edges as for the estimation of the connection probabilities $\gamma_{q\ell}$. {In the same vein, this suggest rates of convergence of the order $n^{-1/2}$ for estimators of the group proportions whereas this should decrease to $n^{-1}$ for estimators of connection probabilities.} {Those two different} rates are observed in a simulation study of \cite{Gazal}. \\
\cite{Latouche_Bayes_EM} and \cite{Come_Latouche} elaborated on this approach in the context of variational Bayes inference, and proposed both a BIC  and an ICL criterion. In this context, the Laplace approximation involved in the classical BIC and ICL criterion is not needed and the corresponding integral can be computed in an exact manner. Note that no formal proof of the consistency of these criteria with respect to the  estimation of $Q$ exist. \\
For SBM, it is most often observed that the difference between ICL and BIC is almost zero. Reminding that this difference corresponds to the conditional entropy of $\bZ$ given $\bY$, this is consistent with the fact that $\pr(\bZ|\bY)$ concentrates around one unique point  \citep{Maria_Matias,Celisse_etal}. \\
More recently, \cite{Channarond} proposed a criterion which {does not} rely on the likelihood or on some approximation of it, but only on the distribution of gaps between the ordered degrees of the nodes. This criterion is proved to be consistent.

\section{Extensions of SBM}
\label{sec:ext}

Several extensions of the SBM have been proposed in the literature with different aims. We present some of them in this section. Default notation are those of the standard SBM defined in \eqref{eq:SBM}.

\subsection{Overlapping groups}
As mentioned earlier, SBM is often used for clustering purposes, that
is to assign individuals (nodes) to groups. As most clustering
methods, the standard SBM assumes that each individual of the
population under study belongs to one single group, as each hidden
state $Z_i$ has a multinomial distribution $\calM(1; \bpi)$ over $\{1,
\dots, Q\}$. This assumption may seem questionable, especially when
analyzing social networks where an individual may play a different
role in each of its relationships with other individuals. Two main
alternatives have been proposed to overcome this limitation of SBMs.

\cite{Airoldi} proposed a mixed-membership model where each node $i$
possesses its own unknown probability vector $\bpi_i = (\pi_{i1},
\dots \pi_{iQ})$. For each link, nodes $i$ and $j$ first choose their
group membership $Z_{i \rightarrow j}$ and $Z_{j \rightarrow i}$ for
this precise link, according to their respective probability vectors
$\bpi_i$ and $\bpi_j$. The value of the link $Y_{ij}$ is then sampled
according to distribution $F(\cdot; \gamma_{Z_{i \rightarrow j}, Z_{j
    \rightarrow i}})$. The number of hidden variables involved in this
model is fairly large (there are $n(Q-1)$ independent $\pi_{iq}$'s and
$2n^2$ independent $Z_{i \rightarrow j}$'s) and a \texttt{MCMC} strategy is
proposed to achieve the inference. A similar mixed-membership model was
proposed by \cite{Erosheva} in the context of a simple mixture model,
without network structure.

\cite{Latouche_overlap} proposed an overlapping version of SBM in
which individuals may belong simultaneously to any subset of classes.
The group membership vector $Z_i = (Z_{iq})_{1 \leq q \leq Q}$ is
drawn as a set of independent Bernoulli variables with respective
probabilities $\pi_q$ (which are not required to sum to 1). Thus
$Z_i$ can take $2^Q$ different values, meaning that one node can belong to zero,
one, two and up to $Q$ classes. In the binary version of the
overlapping SBM, the link $Y_{ij}$ is then present with
logit-probability $Z_i^\intercal W Z_j + Z_i^\intercal U + V^\intercal
Z_j +W^*$, where the matrix $W$, the vectors $U$ and $V$ and the
scalar $W^*$ have to be inferred, as well as the $\pi_q$'s. This
model involves $Q + (Q+1)^2$ parameters, which is much less than its
mixed-membership counterpart. The authors propose a variational
approach for their estimation and for the inference of the membership
vectors $\{Z_i\}_{1\le i \le n}$.

\subsection{Bipartite graphs}
\label{sec:bipartite}

Some networks depict interactions or relationships between two
distinct types of entities, such as authors and journals, chemical
compounds and reactions, hosts species and parasites species, etc. In such networks the link has most often an asymmetric meaning, such as 'published an article in', 'contributes to' or 'is contaminated by'. In such networks, no link between nodes of the same type can exist. When
considering $n$ nodes of the first type and $m$ nodes of the second
type, the adjacency matrix $(Y_{ij})$ is rectangular with $n$
rows and $m$ columns. The SBM model can rephrase in an asymmetric
way, denoting $\{Z_i\}_{1 \leq i \leq n} \in \{1, \dots Q\}^n$ the
memberships of the row nodes and $\{W_j\}_{1 \leq j \leq m} \in \{1,
\dots , K\}^m$  the memberships of the column  nodes. All membership
variables are
drawn independently with multinomial distribution $\calM(1; \bpi)$ for
the $Z_i$'s and $\calM(1; \boldsymbol{\rho})$ for the $W_j$'s. Links $\{Y_{ij}\}$ are
then drawn independently conditional on $\{Z_i, W_j \}_{ 1\le i \le n,
1\le j \le m}$ and with
respective distribution $F(\cdot; \gamma_{Z_i, W_j})$.

This model is actually a latent block model (LBM) first proposed by
\cite{GoN03} in the context of bi-clustering to infer
simultaneously $Q$ row groups and $K$ column groups. The same authors
proposed a variational approximation for the parameter inference. More
recently, \cite{Ker12} proposed a model selection criterion in a variational Bayes context
and \cite{Maria_Matias} proved the convergence of the conditional
distribution of the memberships toward the true ones.

\subsection{Degree-corrected block model}

The regular SBM assumes that the expected connectivity of a node
only depends on the group it belongs to. Indeed, in many situations,
beside their group membership, some nodes may be likely to be more
connected than others because of their individual
specificities. \cite{Karrer} extended the Poisson-valued SBM by taking 
$$
Y_{ij} | Z_i=q, Z_j=\ell \sim \calP(\gamma_{q\ell} \kappa_i \kappa_j),
$$
where $\gamma_{q\ell}$ plays the same role as in the
Poisson-valued SBM and $\kappa_i$ controls the expected degree of node
$i$. A similar model is considered in \cite{Mor09}, who
proposed a relaxation strategy for the inference of the
parameters. \cite{Zhu27092013} proposed an oriented version of this
model and use \texttt{MCMC} for the parameter inference. An asymmetric version
is considered in \cite{RAS11}, generalizing LBM in
the same way.
\cite{YRJ12} proposed  a likelihood-ratio  test for the  comparison of
the   degree  corrected  block   model  with   the  regular   SBM  and
\cite{Zhao_etal} provided  a general characterization  of the community detection criteria
 that will provide a consistent classification under the degree-corrected SBM.

\subsection{Accounting for covariates}

The existence or the value of the links between individuals can sometimes be partially explained by observed covariates. Accounting for such covariates is obviously desirable to better understand the network structure. 
It is first important to distinguish if the covariates are observed at the node level (e.g. $\bx_i=$ (age, sex) of individual $i$) or at the edge level (e.g. $\bx_{ij} = $ (genetic similarity, spatial distance) between species $i$ and $j$). Node covariates $\bx_i = (x_{i1}, \dots, x_{ip})$ are sometimes transformed into edge covariates taking, e.g., $\bx_{ij} = (|x_{i1}-x_{j1}|, \dots, |x_{ip}-x_{jp}|)$. 
\cite{PaR07} propose a brief review of how such an information can be incorporated in a random graph model in absence of hidden structure. For example, in presence of edge covariates, a simple logistic regression model can be fitted as
$$
(Y_{ij})_{i, j} \text{ independent}, \qquad \text{logit}\left(\pr\{Y_{ij} = 1\}\right) = \bx_{ij}^\intercal \bb,
$$
where $\bx_{ij}$ is the vector of covariates and $\bb$ the vector of regression parameters. The statistical inference of this model raises no specific issue.

We now focus on how covariates can be used to enrich SBM. For node covariates $\{\bx_i\}$, \cite{Tallberg05} proposed a multinomial probit model where the distribution of the hidden class $Z_i$ depends on the vector of covariates $\bx_i$
$$
(Z_i)_i \text{ independent}, \quad Z_i \sim \calM(1; \bpi(\bx_i)).
$$
This model states that the covariates act on the edge value through the membership of the nodes. In this context, the author proposed a Bayesian inference approach for which some full conditional distributions can be derived in a close form.

In presence of edge covariates $\{\bx_{ij}\}$, \cite{Mariadassou_10} proposed to combine them with the hidden structure using the generalized linear model framework. In the Poisson case, this leads to 
$$
Y_{ij}  | Z_i=q, Z_j=\ell  \sim \calP(\exp(\bx_{ij}^\intercal  \bb +
\gamma_{q\ell})) .
$$
A similar Gaussian model is considered in \cite{Zanghi_cov}. In both cases, the proposed inference method relies on a variational approximation. 
\cite{Mariadassou_10} argue that the  term related to the hidden state
$\gamma_{q\ell}$ measures the heterogeneity in the network that is not
explained by  the regression term $\bx_{ij}^\intercal  \bb$. A similar
interpretation is given by  \cite{CWA12}, who first apply a regression
step and then perform SBM inference on the residuals.

We mention that our review has not discussed so far the case of \emph{spatial
  networks}, where the nodes correspond to entities that have explicit
geographic locations (a particular case of covariate information).  Integration of this spatial information to
detect  clusters   of  differentially  connected   and  differentially
geographically positioned actors may be of interest in various applications. This has been
done in~\cite{Miele_etal14} in an original way, by considering a SBM with a
regularization term  in the likelihood function that  forces groups to
be the  same for  (geographically) close actors.  This method  is then
applied to modeling of ecological networks.

\subsection{SBM as a graphon model}
\label{sec:graphon}

It   can  be  seen   that  SBM   is  a   graphon  model,   as  defined
through~\eqref{eq:graphon} in Section~\ref{sec:other_latent}, where the function $g$ is block-wise constant, the blocks being rectangular with respective dimension $\pi_q \times \pi_\ell$ and height $\gamma_{q\ell}$. Spliting the interval $[0, 1]$ into $Q$ sub-intervals with respective widths $\pi_1, \dots, \pi_Q$, the hidden state $Z_i$ of SBM is simply the number of the sub-interval into which the corresponding $U_i$ falls.

The inference of the function $g$ has received few attention until now. \cite{Cha12} first proposed a direct estimation of the $p_{ij} = g(U_i, U_j)$ based on a low rank decomposition of the adjacency matrix. More recently, several articles \citep{NIPS2013_5047, 2013arXiv1309.5936W, 2013arXiv1310.6150L,OlW13} proposed to infer $g$ through the parameter of a SBM, which can be seen as a discrete approximation of the graphon. 

It is  important to note that  the graphon model suffers  a strong and
intrinsic identifiability problem as {composing $g$ with} any measure-preserving function from the unit interval to itself will result in the same model. This issue is accounted for in all the papers cited above, but the interpretability of the resulting function is still questionable. As shown in \cite{DiJ07}, subgraphs (also called 'motifs') frequencies are invariant to such transformation and therefore characterize a $W$ graph. The inference of those frequencies is addressed by \cite{2013arXiv1310.6150L}.

\subsection{Network evolution}

More  and more  attention is  paid to  evolutionary networks,  that is
network in which the value of  edges may vary along time. In this this
setting,  the   data  at  hand  consists   in  $\{Y_{ij}(t)\}$  where,
typically, $t$ belongs  to a finite set of  observation times. The set
of nodes is kept fixed here. 

Recently, different proposals have been made to deal
  with groups structure within  dynamic networks. We mention here only
  the  models  connected  with  SBM.  \cite{DBS13}  suppose  that  the
  membership of  the nodes  are kept fixed  and that the  edges values
  evolve according to a conditional  Markov jump process. In the model
  considered by \cite{FSX09} both  the memberships and the conditional
  connection  probabilities  evolve  along time.  \cite{dyn_SBM}  also
  propose a  dynamic version  of SBM where  both node  memberships and
  edges evolve along time. In this paper, the parameter inference relies on a extended Kalman filter algorithm.

\section{Some perspectives}
\label{sec:conclusion}

In this last section, we try to briefly underline what we think are the next
challenges  in the modeling  of heterogeneity  {through latent
  space models} for networks. 

First, it  is important to develop  scalable methods that  are able to
handle  very  large graphs.  Recently,  \cite{Vu_Hunter_S} proposed  a
model for clustering in very large networks: they handled a dataset
with  more than  131,000 nodes  and 17  billions edge  variables. They
consider  discrete-valued  networks,  possibly  with  covariates,  and
assuming independence  of each  dyad (namely $Y_{ij}$  in the  case of
undirected edges and $(Y_{ij},Y_{ji})$ in the case of directed edges),
conditional on the nodes groups.  They rely on the variational
approximation   of   the  \texttt{EM}   algorithm   but  replace   the
\texttt{E}-step by a minorization-maximization algorithm. As a
result,   this  increases   the   lower  bound   in  the   variational
approximation instead of maximizing it at each step. 
{However, the  large network they handle is  still very sparse
  (with about  only 840,798 non zero  edges) and the  number of groups
  they use   $Q=5$ is small with respect to the sample size.}
 Indeed, there is still room for improvement to make model-based clustering methods scalable to very  large  networks.\\

Second, we feel that there is a great need in terms of statistical methods
for the  analysis of  evolutionary networks. Indeed  many technologies
now  give access  to follow-up  observations of  social  or biological
networks. A strong  attention should be paid to  both proper dynamical
modeling and  their associate inference. 
As  shown  throughout  this  paper,  due  to  the  network  structure,
statistical models for statistical  network suffer from very intricate
dependency  structures between  the nodes.  Accounting  for dependency
along time will obviously make it even more complex. The conception of
both  statistically  valid  and  computationally  efficient  inference
methods is an interesting challenge.\\

Finally, statistical properties  of the models and procedures
  should   be   further   studied   from  a   theoretical   point   of
  view.  Asymptotic results have  recently been established but finite
  sample  properties would  also be  welcome.  {Much attention
    has been paid  to the dense case and sparser  setups still need to
    be studied.} 
Validation of the
    procedures  can  not be  limited  to  simulations and  theoretical
    studies  will help better  understand the  models and  thus design
    new inference methods.

\bibliographystyle{abbrvnat}
\bibliography{graphs_review}

\end{document}